\documentclass[11pt]{article}
\usepackage{amsmath,amssymb}
\newcounter{lemma}
\newenvironment{lemma}{\smallskip{\sc Lemma \thelemma.}
\addtocounter{lemma}{1}\em}{\smallskip}
\newenvironment{proof}{\smallskip{\sc Proof.}}{\qed\smallskip}
\newenvironment{generic}[1]{\smallskip{\sc #1.}}{\smallskip}
\newenvironment{genericem}[1]{\smallskip{\sc #1.}\em}{\/\rm \smallskip}
\def\qed{\ifhmode\unskip\nobreak\fi\ifmmode\ifinner\else\hskip.5em\fi\fi
 \hbox{\hskip.5em$\square$\hskip.1em}}
\begin{document}
\title{Approximation of Measures on $S^n$ by discrete Measures}
\author{Marina Nechayeva and Burton Randol}
\date{}
\maketitle

\begin{abstract}We study the asymptotic behavior, as $\rho \rightarrow
\infty$, of discrete measures on $S^{n-1}$ that are induced by
radially projecting point masses concentrated on the integral
lattice-points within dilates $\rho D$ of a compact body $D$, for
various classes of $D$. The results depend sensitively on the
differential geometric properties of $\partial D$.
\end{abstract}

\section{Introduction}

Recently, Douglas, Shiffman, and Zelditch, motivated by
considerations arising from the vacuum selection problem in string
theory, have investigated questions in physics, aspects of which
involve the equidistribution of radial projections of integral
lattice-points on certain hypersurfaces in $R^{n-1}$, e.g., onto
level surfaces of quadratic forms \cite{zelditch1}. One question of
this kind which they take up concerns the study of the asymptotics,
as $\rho \rightarrow \infty$, of the family of discrete measures on
the surface $\partial\mathcal{E}$ of an ellipsoid $\mathcal{E}$
resulting from summing unit-weighted radial projections onto
$\partial\mathcal{E}$ of the lattice-points in $\rho \mathcal{E}-
\{0\}$. If, for example, $\mathcal{E} = S^{n-1}$, this becomes: Does
the action of the measure on $S^{n-1}$ which results from summing
unit-weighted projections onto $S^{n-1}$ of lattice-points $N \in
\Bbb{Z}^n - \{0\}$, with $|N| \leq \rho$, and then rescaling the
result by $n/\rho^n$, tend, on smooth functions, to Lebesgue measure
on $S^{n-1}$ as $\rho \rightarrow \infty$, and if so, how rapidly?
This last question can be regarded as an instance of the general
question of approximating smooth measures on $S^{n-1}$ by discrete
measures, and in this paper, we will examine this question, which
was suggested by questions arising in \cite{zelditch1}, which can be
recast in this form. Depending on the measure being approximated,
issues associated with the curvature of a specific surface
associated with the measure can play a role, in particular the
presence of zones of zero curvature on this associated surface,
which we will illustrate by several examples.

It is a pleasure to mention that our interest in these questions
arose from a conversation with Zelditch about his above-cited recent
work with Douglas and Shiffman.

We will begin the paper with a description of a general method,
which can then be applied to various instances of the problem under
consideration, in particular, to the zero curvature case, which is
not discussed in \cite{zelditch1}. An interested reader could
probably glean the required background from a close reading of one
or more of the papers \cite{zelditch1}, \cite{randol69a},
\cite{randol97}, as well as others, but in our view, expository
clarity makes it very desirable to begin with such an exposition,
which we have written to lead as smoothly as possible into the
notation and approaches of the various papers to which we refer. In
particular, since detailed treatments of instances of the zero
curvature cases are in general quite prolonged, and very closely
mimic existing discussions in the literature for the corresponding
classical constant-density lattice-point problem, we have confined
our discussion of these examples to general descriptions of how the
overall techniques described in this paper can be adapted to closely
follow treatments of the corresponding constant-density cases in the
literature.

\section{}
\label{sec_A}

Suppose $d\mu$ is a Borel measure on $S^{n-1}$ having a continuous,
positive, piecewise smooth density function $m$ with respect to
Lebesgue measure. If $d\mu$, as in the example above, happens to be
Lebesgue measure, i.e., $m = m(\theta) \equiv 1$, the corresponding
family of discrete measures is parametrized by $\rho > 0$, and
supported, as indicated above, on the unit-weighted radial
projections of the non-zero integral lattice-points $N$, for which
\mbox{$|N| < \rho$}. Neglecting rescaling, the natural counterpart
for a general measure on $S^{n-1}$ having positive density $m$ with
respect to Lebesgue measure is the family $d\Gamma_\rho$ of discrete
measures supported on the unit-weighted radial projections of the
non-zero integral lattice-points in $\rho D$, where $D$ is the
compact set whose boundary is given in polar coordinates by $r =
(m(\theta))^{1/n}$. We assume that $m(\theta)$, or equivalently
$\partial D$, is sufficiently regular so that the divergence theorem
is valid for $D$ -- for example, $\partial D$ could be smooth, or
the boundary of a polyhedron.

\section{}

We begin by noting that the effect of the above projection measure
on a smooth function $f(\theta)$ on $S^{n-1}$ is identical to that
of the lattice-point count over $\rho D -\{0\}$, weighted by the
weight-zero homogeneous extension $F$ of $f$ to $R^n -\{0\}$. There
are various analytical approaches to the estimation of such sums,
e.g., Riesz means, as in  \cite{hlawka1}, \cite{randol66a}, etc., or
convolution smoothing techniques, as in, for example,
\cite{zelditch1} and \cite{randol69a}. We will employ a convolution
smoothing technique, since it is relatively simple to implement and
describe. In our outline of this approach, we will temporarily
assume, for minor technical reasons, that $D$ has been modified by
the removal of a neighborhood of the origin. Then, assuming for the
moment that $f$ is positive, and denoting its homogeneous extension
to $R^n$ by $F$, the above lattice-point count $N_D(\rho)$ in $\rho
D$ can be overestimated by summing over the weighted lattice-points
in a slight expansion $E(\rho)$ of $\rho D$ and underestimated by
summing over the weighted lattice-points in a slight contraction
$C(\rho)$ of $\rho D$. The Poisson summation formula provides a
natural analytic approach for handling such sums, but since the
function $F$ is not smooth when multiplied by the indicator
functions of $E(\rho)$ or $C(\rho)$, i.e., at the boundaries of
$E(\rho)$ and $C(\rho)$, we will employ a slight refinement of the
above idea, which can be approximately described by saying that we
convolve $F$, restricted to $E(\rho)$ and $C(\rho)$, respectively,
with a smooth compactly supported approximate delta function
$\delta_\epsilon$, whose support is chosen to lie within a ball of
radius $\epsilon (\rho)$ which is sufficiently small so that the
convolutions are very close to the function $F$ restricted to $\rho
D$, except in small neighborhoods of the boundary. If the support of
$\delta_\epsilon$ is of sufficiently small diameter, it will then be
exactly true, if $F$ is constant, and approximately true otherwise,
that $\delta_\epsilon * F_{C(\rho)} \leq F_{\rho D} \leq
\delta_\epsilon * F_{E(\rho)}$, where $F_{\mathfrak{S}}$ denotes the
product of $F$ with the indicator function $\chi_\mathfrak{S}$ of
$\mathfrak{S}$. If $F$ is not constant, the deviation from
correctness of the above inequality can be quantified in terms of a
bound for the directional derivatives of $F$ and the diameter of the
support of $\delta$. In view of this, if, on $E(\rho)$ and $C(\rho)$
respectively, we modify $F_{E(\rho)}$ and $F_{C(\rho)}$ by the
addition and subtraction of a suitably small quantity whose size
depends on $\rho$, the above inequality becomes correct, and can be
exploited, following the method of \cite{randol69a}, which treated
the case of constant density, to derive estimates for $N_D(\rho)$.
If $f$ is of mixed sign, we can express it as the difference of two
positive smooth functions, proceed as above, and later combine the
estimates for both. It is almost immediate that any smooth function
is expressible as the difference of two positive smooth functions,
whose size and the size of whose derivatives are comparable to those
of the original function. For instance, we could define $f = f^+ -
f^-$, where $f^+ = f + c$,
 and $f^- = c$, where $c >|\mathrm{min}_{\theta \in S^{n-1}}
 f(\theta)|$.

We now show in greater detail how to carry out this program. The
outcome in a general sense will be that after rescaling by
$n/\rho^n$, discrete sphere measures which arise in the above way
will converge to $d\mu$ on $S^{n-1}$ at the same rate as that at
which $1/\rho^n$ times the standard unit-density lattice-point sum
for $\rho D$ converges to the measure of the set $D$. The specific
sense in which this is the case will be developed more precisely
below.

\section{}

We will retain the notation of Section (\ref{sec_A}). As remarked
there, the task is equivalent to the estimation, as $\rho
\rightarrow \infty$, of the lattice-point count over $\rho D$, with
lattice-points counted with density $F$, where $F$ is a smooth
homogeneous function of weight 0 on $R^n - \{0\}$ (if $F$ is
constant, it is smooth at the origin as well). Unless $F$ is
constant, the origin causes a small, easily surmounted technical
difficulty, which can be handled in various ways. We will, following
Douglas, Shiffman, and Zelditch, deal with this by estimating the
weighted lattice-point count in the shell $\rho S = \rho D -
\frac{\rho}{2}D$, which is a dilate of the basic shell $S = D -
\frac{1}{2} D$. Since this estimate will be uniform in the dilation
parameter, this automatically leads to estimates for the count for
$\frac{1}{2} \rho S$, $\frac{1}{4} \rho S$, $\frac{1}{8} \rho
S$,$\dots\,$. We then obtain the desired result for $\rho D$ by
adding up the results for $\rho S, \frac{1}{2} \rho S, \frac{1}{4}
\rho S, \frac{1}{8} \rho S, \dots\,$. The weighted lattice-point
count for $\rho S$ will be estimated along the lines outlined above,
i.e., by bracketing the true count between appropriate convolutions.

\section{}

We now pass to the details.

We will use the following lemma (cf.\ \cite{randol66a}, pp.\
262--263), whose role is to show that the Fourier transform
$\hat{g}(y)$ of a smooth function $g(x)$ on a domain $\mathfrak{S}$
for which the divergence theorem is valid can be expressed as the
Fourier transform of a smooth function on $\partial \mathfrak{S}$,
and that the process of transfer to the boundary picks up a factor
of $1/|y|$.

\begin{lemma} Suppose $g(x)$ is a smooth function on $R^n$, and $y$ a fixed vector in
$R^n$. Then there exists a smooth vector field $\mathbf{F(x)}$ on
$R^n$, such that

\begin{enumerate}

\item $\mathbf{div}\;[(2\pi i|y|)^{-1} e^{2\pi i(x,y)}\mathbf{F}(x)] =
e^{2\pi i(x,y)} g(x)$

\item The derivatives up to order $k$ of the components of $\mathbf{F}(x)$
can be bounded, independently of $y$, in terms of bounds for the
derivatives up to order $k+1$ of $g(x)$.

\end{enumerate}

\end{lemma}

\begin{proof} \begin{sloppypar} This is proved in Lemma 3 of \cite{randol66a}.
In particular, since, setting $x=(x_1,\dots ,x_n)$, $y=(y_1, \ldots
, y_n)$, and $\mathbf{F}=(F_1,\ldots ,F_n)$, the divergence of
$(2\pi i|y|)^{-1} e^{2\pi i(x,y)}\mathbf{F}(x)$ is
\begin{eqnarray}\lefteqn{ e^{2\pi i(x,y)} [((y_1/|y|)F_1 + (2\pi
i|y|)^{-1}\partial F_1/\partial x_1) + \ldots}\nonumber\\ & & +
((y_n/|y|)F_n + (2\pi i|y|)^{-1}\partial F_n/\partial x_n)]
\,,\nonumber\end{eqnarray}\end{sloppypar}

\[ = e^{2\pi i(x,y)}[(y/|y|\,,\mathbf{F}) + (2\pi i|y|)^{-1}\mathbf{div \, F}] \,,\] the requirements
of the present lemma will be satisfied if

\begin{enumerate}

\item $\mathbf{div \, F} = 0$

\item $(y/|y|\,,\mathbf{F}) = g$

\end{enumerate}with the derivatives of the $F_j$'s bounded as indicated
above. This is, however, precisely the assertion of Lemma 3 of
\cite{randol66a}, taking $\beta = y/|y|$ in the notation of that
lemma.\end{proof}

\begin{generic}{Corollary} The first assertion of the lemma immediately
implies, by the divergence theorem, that \[ \int_\mathfrak{S}
e^{2\pi i (x,y)} g(x)\, dx = (2\pi i|y|)^{-1} \int_{\partial
\mathfrak{S}} e^{2\pi i (x,y)} (\mathbf{F}(x),n(x))\, ds_x \,,\]
where $n(x)$ is the exterior normal to $\partial \mathfrak{S}$ at
$x$. This gives the desired expression of the Fourier transform over
$\mathfrak{S}$ in terms of a Fourier transform over $\partial
\mathfrak{S}$.\end{generic}

As previously indicated, we will assume that $F>0$, since the
general case can be deduced from this one by taking differences of
the resulting estimates. For a small positive parameter
$\epsilon(\rho)$, which will depend on $\rho$, and which we will
sometimes simply write as $\epsilon$ for short, we define,
respectively, a slight expansion $E(\rho)$ of $\rho S$, and a slight
contraction $C(\rho)$ of $\rho S$, by setting $E(\rho) =
\mbox{$(\rho +\epsilon)D - (\frac{\rho}{2} -\epsilon)D$}$, and
$C(\rho) = (\rho - \epsilon)D - (\frac{\rho}{2} + \epsilon)D$. Since
we are interested in the asymptotics of the weighted lattice-point
count over $\rho S$ as $\rho \rightarrow \infty$, we may assume, by
replacing $S$ by a sufficiently large dilate of its original self if
necessary, that for small $\epsilon (\rho)$, $E(\rho)$ contains an
$\epsilon (\rho)$-neighborhood of $\rho S$, and that $\rho S$
contains an $\epsilon (\rho)$-neighborhood of $C(\rho)$.

Throughout the following, our basic approach is via the Poisson
summation formula, and is essentially the same technique that was
employed in \cite{randol69a}, the only difference being that now the
weight function is no longer a constant. Since $N_S(\rho) =
\sum_{N\in \rho S} F(N)$, a purely formal application of the Poisson
summation formula leads to
\begin{equation} \rho ^n \sum_{N\in \mathbb{Z}^n} {\hat{F}}_S (\rho
N) \label{Poisson1}\end{equation} as an analytic expression for the
weighted lattice-point count, where ${\hat{F}}_\mathfrak{S}$ denotes
the Fourier transform of $F_\mathfrak{S}$.

As previously mentioned, there are convergence difficulties, since
the function which is equal to $F(x)$ on $S$ and equal to zero on
the complement of $S$ is discontinuous at the boundary of $S$ unless
$F \equiv 0$. This problem is of course present in the constant
density case as well, and it can be addressed, as was done in
\cite{randol69a} for that case, and as we have previously indicated,
by slightly expanding and contracting $\rho S$, into $E(\rho)$ and
$C(\rho)$, respectively, and then smoothing the indicator functions
of $E(\rho)$ and $C(\rho)$ by convolving them with a smooth
approximate delta function $\delta_\epsilon (x)$, whose support lies
in a ball of diameter $\epsilon(\rho)$. Then if $F$ is constant,
\begin{equation} \sum_{N\in C(\rho)}(\delta_\epsilon * F_C)(N) \leq
N_S(\rho) \leq \sum_{N\in E(\rho)}(\delta_\epsilon * F_E)(N)
\,,\label{ineq1}\end{equation} and since the functions being summed
are smooth with compact support, the Poisson summation formula can
be applied to both sides of the above inequality to obtain
asymptotics for $N_S(\rho)$ (cf.\ \cite{randol69a}).

In the case of non-constant $F(x)$, a minor technical issue
connected with this approach arises from the fact that $F(x)$ is no
longer necessarily equal to its convolution by $\delta_\epsilon
(x)$. However, for large $\rho$, $(\delta_\epsilon * F)(x)$ is very
close to $F$ on a dilated shell $\rho \mathfrak{S}$, since any first
derivative of $F$ is homogeneous of weight $-1$, and therefore the
oscillation of $F$ over the intersection of a ball of diameter
$\epsilon$ with $\rho\mathfrak{S}$ is at most of the order of
$(1/\rho)\epsilon(\rho)$. I.e., $F$ is close to being constant on
$\rho \mathfrak{S}$ as $\rho$ becomes large, so, as previously
indicated, the addition and subtraction of a suitable small quantity
$\eta_\rho$ in inequality (\ref{ineq1}) for the constant density
case results in a correct inequality for the general case:

\[ \sum_{N\in C(\rho)}[\delta_\epsilon * (F_C - \eta_\rho\chi_C)](N) \leq
N_S(\rho) \leq \sum_{N\in E(\rho)}[\delta_\epsilon * (F_E
+\eta_\rho\chi_E)](N) \,,\] or

\begin{eqnarray} \lefteqn{\sum_{N\in C(\rho)}(\delta_\epsilon * F_C)(N) -
\sum_{N\in C(\rho)}(\delta_\epsilon * \eta_\rho\chi_C)(N)  \leq N_S(\rho)}\nonumber\\
& & \leq \sum_{N\in E(\rho)}(\delta_\epsilon * F_E)(N) + \sum_{N\in
E(\rho)}(\delta_\epsilon *
\eta_\rho\chi_E)(N)\,,\label{ineq2}\end{eqnarray} where
$\chi_\mathfrak{S}$ denotes the indicator function of
$\mathfrak{S}$, and for notational simplicity we have written $C$
and $E$ for $C(\rho)$ and $E(\rho)$, respectively.

Since $\eta_\rho$ will be chosen to be of the order of
$(1/\rho)\epsilon(\rho)$, the analysis of the contribution to the
above inequality that is attributable to the presence of the
$\eta_\rho$ terms on either side will be easily seen to lead to the
presence, on either side, of a term of the order of
$\rho^{n-1}\epsilon(\rho)$. The choice of $\epsilon$ as a function
of $\rho$ will depend on the shape of $D$, and for the time being,
we will simply carry it along as an unspecified function of $\rho$.

When we come to the analysis of the sums in (\ref{ineq2}), the heart
of the matter lies in the asymptotics of the Fourier transform of
$F$ over $C(\rho)$ and over $E(\rho)$, which reduces, by an obvious
rescaling, to the analysis of the asymptotics of the Fourier
transform ${\hat{F}}_{\frak{S}}$ of $F$ over shells $\frak{S}$ which
are very close to $S$ for large $\rho$.

Accordingly, we next take up the asymptotic analysis of \[
{\hat{F}}_{\frak{S}}(y) = \int_{\frak{S}} F(x) e^{2\pi i(x,y)}\,dx
\,,\] where $\frak{S} = \frak{S}(\rho)$ can stand for a rescaling by
$1/\rho$ of either $C(\rho)$ or $E(\rho)$.

By Lemma 1, \begin{equation} \int_{\frak{S}} e^{2\pi i (x,y)} F(x)\,
dx = (2\pi i|y|)^{-1} \int_{\partial \frak{S}} e^{2\pi i (x,y)}
(\mathbf{F}(x),n(x))\, ds_x \,,\label{eq1}\end{equation} where the
size of the derivatives of the vector field $\mathbf{F}$ are
controlled by the size of the derivatives of the homogeneous
function $F$, which, since $\mathfrak{S}$ is  a shell, we can smooth
at the origin without affecting its values in $\mathfrak{S}$.

The boundary of $\frak{S}$ splits into two components, which are
analyzed similarly. Thus, the asymptotic analysis of \[
\int_{\frak{S}} F(x) e^{2\pi i(x,y)}\,dx \,\] amounts to the study
of the asymptotic behavior of integrals \[ \int_{\mathcal{S}} g(x)
e^{2\pi i(x,y)}\,dx \,,\] where $\mathcal{S}$ is a hypersurface in
$R^n$ of the considered type, and $g$ is a smooth function on
$\mathcal{S}$. This is a highly studied subject (cf.\
\cite{hlawka1}, \cite{randol69a}, \cite{randol69b}, et al.). The
results are very dependent on the geometry of $\mathcal{S}$,
although the case in which $\mathcal{S}$ has everywhere positive
Gaussian curvature is easy to analyze, using, for example, a
straightforward application of stationary phase, or even a simple
calculus argument. It is well-known that in this case (cf.\
\cite{hlawka1}), \begin{equation} \int_{\mathcal{S}} g(x) e^{2\pi
i(x,y)}\,dx \ll |y|^{-(n-1)/2}\,,\label{FT_ineq}\end{equation} where
the implied constant depends on bounds for the derivatives of $g$,
which, for the purposes of this discussion, can be thought of as a
function defined in a neighborhood of $\mathcal{S}$. By the
corollary to Lemma 1, it follows that in the positive curvature
case, for an $n$-dimensional shell $\mathfrak{S}$, \begin{equation}
\int_{\frak{S}} e^{2\pi i (x,y)} F(x)\, dx \ll
|y|^{-(n+1)/2}\,.\label{FT_ineq2}\end{equation} Although the
positive curvature case has already been discussed in
\cite{zelditch1}, we will use this case to illustrate our approach,
since it provides an exceptionally simple paradigmatic description
of techniques of considerably more general applicability. We will
then describe modifications that are required in illustrative cases
for which $\partial D$ does not have everywhere positive curvature.
In general, specific applications require individual adaptations of
the method, but these adaptations can often be easily inferred from
existing treatments of the constant-density case.

The technique closely follows that of \cite{randol69a}, and so will
be simply outlined here. Denoting as above by $N_S(\rho)$ the
weighted lattice-point count in $\rho S$, we will indicate how to
use the Fourier transform bound to estimate the right side of the
inequality (\ref{ineq2}). The left side is handled in a similar way.
The Poisson summation formula can be applied to the right side of
(\ref{ineq2}). We first take up the term \[\sum_{N\in
E(\rho)}(\delta_\epsilon * F_E)(N)\,,\] which, by the Poisson
summation formula equals
\[\sum_{N\in E(\rho)} {\hat{\delta}}_\epsilon (N)\,
{\hat{F}}_{E(\rho)}(N)\,,\]
\[\;\;\; = {\hat{\delta}}(0)\,\hat{F}(0)  +
\sideset{}{'}\sum_{N\in E(\rho)}{\hat{\delta}}_\epsilon (N)\,
{\hat{F}}_{E(\rho)}(N)\,,\] where the prime on the summation sign
means that the origin is omitted from the sum.

Since $\hat{\delta}(0) = 1$, the last quantity can be rewritten as
\[ \int_{E(\rho)} F(x)\,dx  + \sideset{}{'}\sum_{N\in
E(\rho)}{\hat{\delta}}_\epsilon (N)\, {\hat{F}}_{E(\rho)}(N)\,.\]

If in the above we replace the term $\int_{E(\rho)} F(x)\,dx$ by
$\int_{\rho S} F(x)\,dx$, the error corresponding to the replacement
can be estimated by noting that \begin{eqnarray*}
\left|\int_{E(\rho)}
F(x)\,dx - \int_{\rho S} F(x)\,dx\right| & = & \left|\int_{E(\rho) -\rho S}F(x)\,dx\right|
\\ \\ & \ll & \mathrm{vol}(E(\rho) -\rho S)\\
\\ & \ll & \epsilon\rho^{n-1}\,.\end{eqnarray*}

The remaining term on the right side of (\ref{ineq2}), namely
\[\sum_{N\in E(\rho)}(\delta_\epsilon * \eta_\rho\chi_E)(N)\,,\] is clearly $\ll$
($\eta_\rho$ times the volume of $E_\rho$), i.e., of order
$\eta_\rho\, \rho^n$.

 We can thus rewrite the right side of (\ref{ineq2}) as \begin{equation}
 \int_{\rho S} F(x)\,dx  +
O(\epsilon\rho^{n-1} + \eta\rho^n) + \sideset{}{'}\sum_{N\in
E(\rho)}{\hat{\delta}}_\epsilon (N)\,
{\hat{F}}_{E(\rho)}(N)\,.\label{expression1}\end{equation}

The last sum is handled in exactly the same way as in the well-known
case in which $F$ is constant. Namely, the necessary estimate on
$\hat{F}_{E(\rho)}(N)$ is provided by (\ref{FT_ineq2}), and is
identical with the estimate for the constant case. As in that case,
integration by parts shows that the term ${\hat{\delta}}_\epsilon
(N)$ is $\ll 1/(1+\epsilon |N|)^k$, for any fixed integer $k$. The
sum is estimated by splitting it into two parts, over $N$ for which
$|N| < 1/\epsilon$, and over $N$ for which $|N| \geq 1/\epsilon$,
respectively. Each of these parts is estimated by comparison with an
integral, using a sufficiently high value of $k$ to produce
convergence, and the estimate (\ref{FT_ineq2}) for the Fourier
transform. The result is then minimized by choosing $\epsilon$ to
balance the estimates, which leads to the choice of $\epsilon =
\rho^{-(n-1)/(n+1)}$, which in turn results in an estimate of
$\rho^{(n-1)(n/(n+1))}$ for the infinite sum in (\ref{expression1}).
The above choice of $\epsilon$ results in the same estimate for the
term of order $\epsilon\rho^{n-1}$ in (\ref{expression1}), and since
$\eta$ is of the order of $\rho^{-1}\epsilon$, we obtain the same
estimate for the term $\eta\rho^n$ in (\ref{expression1}).

Since a similar analysis can be applied to the left side of the
inequality (\ref{ineq2}), this shows that in the positive curvature
case, the weighted lattice-point count over $\rho S$ equals
\[\int_{\rho S} F(x)\,dx + O(\rho^{(n-1)(n/(n+1))})\,.\]

In order to analyze the weighted lattice-point count over $\rho D$,
we add up repetitions of the above quantity, corresponding to $\rho,
\rho /2, \rho /2^2, \ldots ,\rho /2^k$, where $k$ is taken to be of
the order of $\log_2 \rho$, so that $\rho /2^k < 1$. This takes
account of all lattice-points in $\rho D$ except perhaps for a fixed
finite number near the origin, and since these latter have no effect
on the asymptotics, we find that \begin{eqnarray*} N_D(\rho) & =
 & \int_{\rho D} F(x)\,dx + \sum_{j=0}^k O((\rho
/2^j)^{(n-1)(n/(n+1))})\\\\ & = & \rho^n \int_D F(x)\,dx +
\sum_{j=0}^k O((\rho /2^j)^{(n-1)(n/(n+1))})\,.\end{eqnarray*} Since
the implied constants in the $O$ terms are uniformly bounded, this
implies that \[ N_D(\rho)  = \rho^n\int_{D} F(x)\,dx +
O(\rho^{(n-1)(n/(n+1))})\,.\]

Now it is a consequence of equation (4) of \cite{randol66a} that
\begin{equation} \int_{D} F(x)\,dx = (1/n)\int_{S^{n-1}}
f(\theta)m(\theta)\,d\theta \label{expression2}\,,\end{equation} so
we obtain the following equivalent form of the result of Douglas,
Shiffman, and Zelditch in the positive curvature case:

\bigskip

\begin{genericem}{Theorem} {\rm (cf.\ \cite{zelditch1})} With notation as above,
 \[\int_{S^{n-1}} f(\theta)m(\theta)\,d\theta  =
(n/\rho^n) N_D(\rho) + O(\rho^{-2n/(n+1)})\,.\]\end{genericem}

\noindent This is an estimate for the rapidity of the convergence of
the discrete measures $(n/\rho^n)\,d\Gamma_\rho$ to \ $d\mu =
m(\theta)d\theta$. The estimate, which was given in a somewhat
different form in \cite{zelditch1}, coincides with the usual scaled
error term for the standard lattice-point problem for bodies with
boundaries having strictly positive Gaussian curvature (cf.\
\cite{hlawka1}).

To justify (\ref{expression2}), we note that from (4) of
\cite{randol66a} it follows that \[ \int_D F(x)\,dx = \int_0^1
t^{n-1}\,dt\,\int_{\partial D} F(tx) (x,n(x))\,ds_x\,,\] where
$n(x)$ is the outward normal to $\partial D$. Since $F$ is
homogeneous of weight $0$, the double integral equals
\[(1/n)\int_{\partial D} F(x)\,(x,n(x))\,ds_x \,,\] and
coordinatizing $\partial D$ by $S^{n-1}$ via the radial map, this
becomes \[ \int_{S^{n-1}}
f(\theta)\,\Phi(\theta)\,(x,n(x))\,d\theta\,,\] where $\Phi(\theta)$
is the Radon-Nikodym derivative $ds_x/d\theta$, and $x$ is regarded
as a function of $\theta$ via the radial map. It follows easily from
elementary geometric considerations that
\[\Phi(\theta) = |x|^n/(x,n(x))\,,\] so
\begin{eqnarray*}\int_{\partial D} F(x)\,(x,n(x))\,ds_x  & = &
\int_{S^{n-1}}f(\theta)\,\Phi(\theta)\,(x,n(x))\,d\theta\\\\ & = &
\int_{S^{n-1}} |x|^n\,f(\theta)\,d\theta\,. \end{eqnarray*}

But on $S^{n-1}$, $|x| = (m(\theta))^{1/n}$, so the last integral
equals \[\int_{S^{n-1}} f(\theta)\,m(\theta)\,d\theta\,,\] from
which (\ref{expression2}) follows.

\section{The Case in which Curvature can Vanish}

Classical lattice-point asymptotics for dilates $\rho D$ of a body
$D$ for which the curvature of $\partial D$ is not always positive
can be quite intricate, and are highly dependent on the manner in
which the curvature of $\partial D$ vanishes, as well as on the
placement of $D$ in relation to the integer lattice. For example, in
$\mathbf{R}^2$, if the curvature of $\partial D$ vanishes to finite
order at a finite number of points, the classical lattice-point
asymptotics depend on the order to which the curvature vanishes at
the points in question, as well as on Diophantine properties of the
normal vectors to $\partial D$ at those points (cf.\
\cite{randol69a}). The situation can become much more complicated in
higher dimensions, and nuances of this type similarly affect
convergence rates of the discrete measures with which we are dealing
in this paper. Since, once one has the necessary Fourier transform
asymptotics, methods for dealing with these issues closely mimic
those for the classical case, we will content ourselves with briefly
indicating what happens in a few interesting representative cases.

As indicated above, the central analytic issue is the detailed
asymptotics of the Fourier transform of smooth functions on
$\partial D$. In the case of everywhere positive curvature, one has
the previously mentioned result that \[ \int_{\mathcal{\partial D}}
g(x) e^{2\pi i(x,y)}\,dx  \ll |y|^{-(n-1)/2}\,,\] where the estimate
does not depend on the directional component of $y$.

This estimate is generally false if the curvature of $\partial D$
vanishes on some non-void subset of $\partial D$, and a useful
description of what happens can be quite complicated, for example,
if \[\partial D = \{(x_1,\ldots ,x_n)\;|\; {x_1}^{2k} + \cdots +
{x_n}^{2k} =1\}\,,\] or if $D$ is a polyhedron, (cf.\
\cite{randol66b} and \cite{randol84b}).

If we write \[\int_{D} e^{2\pi i(x,y)}\,dx\] in polar coordinates as
$\Psi(r,\phi)$, one quite general fact along these lines is given by
Theorem 1 of \cite{randol69b}. Namely, if $\partial D$ is
real-analytic and $D$ is convex, then the function
\[\Lambda(\phi) = \sup_r r^{(n+1)/2}\, \Psi(r,\phi)\] is in
$\mathbf{L}^p (S^{n-1})$, for some $p >2$. Thus, under conditions of
considerable generality, the Fourier transform asymptotics coincide
with those of the positive curvature case, up to multiplication by a
function in $\mathbf{L}^p (S^{n-1})$. The convexity hypothesis is
unnecessary in $\mathbf{R}^2$, and possibly in higher dimensions,
although this is not generally known. The requirement of
real-analyticity can be replaced by  somewhat weaker hypotheses
(cf.\ \cite{svensson}). For later related papers, cf.\
\cite{brandolini1}, \cite{varchenko1}.

\section*{Example 1} The arguments by which the above-mentioned theorem is proved
can be uneventfully applied to obtain the following counterpart of
(\ref{FT_ineq}), in the case, for example, in which $D$ is convex
and $\partial D$ is real-analytic (and under weaker hypotheses in
$R^2$): \begin{equation} \int_{\partial D} g(x)e^{2\pi i (x,y)}\, dx
\ll \Lambda(\phi)
|y|^{-(n-1)/2}\,,\label{L_p_estimate}\end{equation} where $\Lambda
(\phi)$ is in $\mathbf{L}^p (S^{n-1})$, for some $p
>2$. This implies, by a straightforward application of the
techniques of this paper, that if we modify the definition of the
discrete measures $\Gamma_\rho$ by replacing $\mathbf{Z}^n$ by its
image under the action of an \mbox{element $\gamma$} of
$\mathbf{SO}(n)$, call the resulting measures $\Gamma_\rho(\gamma)$,
and denote the sum corresponding to $N_D(\rho)$ by
$N_D(\rho,\gamma)$, that we obtain the following result:

\bigskip

\begin{genericem}{Theorem A (vanishing curvature case)} With notation
as above,\[ \int_{\mathbf{SO}(n)} |R_D(\rho,\gamma)|\,d\gamma \ll
O(\rho^{-2n/(n+1)})\,,\] where
\[R_D(\rho,\gamma) = \int_{S^{n-1}} f(\theta)\,m(\theta)\,d\theta - (n/\rho^n)N_D(\rho,\gamma)\,.
\]\end{genericem}

I.e., even for a quite general version of the zero curvature case,
the error estimate for the positive curvature case holds for the
$L^1$ norm over the rotation group of the errors for the ``rotated''
measures. There are, of course, various consequences corresponding
to other $L^p$ norms as well.

\section*{Example 2} Our next example is the special case in which \[\partial D =
\{(x_1,\ldots ,x_n)\;|\; {x_1}^{2k} + \cdots + {x_n}^{2k} =1\}\,.\]
The relevant estimate on the Fourier transform is given by Theorem 2
of \cite{randol66b}, which states that for sufficiently smooth $g$,
if we express
\[\int_{\mathcal{\partial D}} g(x) e^{2\pi i(x,y)}\,dx\] in polar
coordinates as $\Psi(r,\phi)$ ($\phi \in S^{n-1}$, $\phi =
(\phi_1^*,\ldots ,\phi_n^*)$), then on the set of points $(r,\phi)$
for which exactly $j$ of the $\phi_i^*$'s vanish,
\begin{equation}\Psi(r,\phi) \ll (A(\phi))^{-\beta}\,r^{-\alpha_j}\,,\label{zero_curvature_est}
\end{equation} where in the above, $A(\phi)$ is the product
of the non-zero $\phi^*_j$'s, $\beta = (k-1)/(2k-1)$, and $\alpha_j
= (j/2k) + (n-j-1)/2$. This is, of course, in some sense a special
case of (\ref{L_p_estimate}), but it gives considerably more
detailed information about the asymptotics of the Fourier transform
for this case.

The classical lattice-point problem for this case is discussed in
\cite{randol66b}. The principal result is that the error term is of
order $\rho^R$, where $R = \max(A,B)$ with $A= (2k-1)(n-1)/2k$ and
$B= n(n-1)/(n+1)$.  The estimate is best possible if $A
>B$. This result is obtained by separately analyzing groups of
lattice-points on the Fourier transform side of the Poisson
summation formula, where the grouping is arranged into sets defined
by the vanishing of a specific number of the $\phi_i^*$'s. That the
result is sometimes optimal is a consequence of the observation that
if $A >B$, the contribution to the Fourier transform side of the
Poisson summation formula coming from lattice-points on axes defined
by the vanishing of all but one of the $\phi_i^*$'s constitutes the
major contribution to the error, and that the behavior of this
contribution can be analyzed by a stationary phase argument.

The adaptation of this result to the present context is
straightforward, and the arguments are nearly identical to those for
the original result, so as before, we will content ourselves with a
statement of the theorem in the present context, using notation as
above. The result in the present context is:

\bigskip

\begin{genericem}{Theorem B (zero curvature case)} With notation as above, if \[\partial D =
\{(x_1,\ldots ,x_n)\;|\; {x_1}^{2k} + \cdots
+ {x_n}^{2k} =1\}\,,\] then
\[\int_{S^{n-1}}f(\theta)m(\theta)\,d\theta  = (n/\rho^n) N_D(\rho) +
O(\rho^{A-n})\,,\] and if $A>B$, this is best
possible.\end{genericem}

\bigskip

\section*{Example 3} As a final example, we will consider the polyhedral case, and
in the interests of expository and combinatorial simplicity,
describe a typical result for the 2-dimensional case.

In the classical constant-density lattice-point problem, if at least
one of the perpendicular vectors to a face of a compact polyhedron
has rational coordinates, there are an infinite number of
\mbox{$\rho_i \rightarrow \infty$} from which an infinitesimal
displacement results in a modification of the lattice-point count of
order $\rho^{n-1}$, so in this circumstance the error estimate is of
true order $\rho^{n-1}$. Since a simple estimate of Gauss shows that
the error term is always $\ll \rho^{n-1}$, polyhedra can be worst
possible cases for lattice-point error asymptotics.

\begin{sloppypar}Paradoxically, this situation is not generic for polyhedra, as
was noticed by Khintchine \cite{khintchine} in the 2-dimensional
case. Khintchine's result is that the error estimate corresponding
to almost any rotation of the integer lattice $\mathbb{Z}^2$ is in
fact extremely small. Specifically, for any $\epsilon >0$, it is
almost always $\ll \log^{1+\epsilon}\rho$. Later work
\cite{randol69a}, \cite{randol84b}, \cite{randol97},
\cite{skriganov93}, \cite{skriganov98}, \cite{skriganov00},
\cite{tarnopol2} has considered various aspects of the
$n$-dimensional case, as well as additional features and refinements
of the 2-dimensional case.\end{sloppypar}

As in the previous examples, methods for the classical
constant-density lattice-point problem carry over with very little
change to our context. As a typical illustrative representative of
what can be expected, we begin by recalling the result, due to
Skigranov \cite{skriganov93}, that if $D$ is an algebraic polygon,
(one for which the ratios of the direction numbers of the normals to
its faces are all algebraic of degree $\geq 2$), then the classical
lattice-point error term is $\ll \rho^\epsilon$, for any $\epsilon
>0$. This result is also a special case of theorems in later papers
\cite{randol97} and \cite{skriganov98}. It can be obtained by using
a more detailed form of the estimate (\ref{L_p_estimate}), as was
done in the previous example (\ref{zero_curvature_est}), but one
which is adapted to the particular case in which $D$ is a polygon
(cf.\ \cite{randol66a}). With such an estimate in hand, the
lattice-points on the Fourier transform side of the Poisson
summation formula are split into two groups: those in finite-width
bands surrounding the normal vectors to the sides of $D$, and all
the rest. The contribution from the lattice-points exterior to the
bands can be estimated by comparison with an integral, while the
series arising from the contributions from lattice-points within
bands is estimated by using Diophantine properties of the ratios of
direction numbers associated to the corresponding normals. Since the
relevant estimate for the Fourier transform of $D$ is singular at
these directions, the poor approximability of these ratios, which is
a consequence of Roth's Theorem, is crucial (cf.\ e.g.,\
\cite{randol84b}, p.\ 858 for a similar argument). The corresponding
result, in the context of the present paper and expressed in the
notation of this paper, is that for a polygon of the above
type,

\bigskip

\begin{genericem}{Theorem C (polygonal case)} With notation as above,
\[\int_{S^{n-1}}f(\theta)m(\theta)\,d\theta =
(n/\rho^n)N_D(\rho) + O(\rho^{\epsilon -n})\,.\]\end{genericem}

\section{Conclusion}

We have described a general method for describing the accuracy with
which a large class of measures on $S^n$ can be approximated by a
naturally associated family of discrete measures. The case in which
$\partial D$, in the notation of this paper, has everywhere positive
curvature has been previously studied in \cite{zelditch1}, and is
taken in this paper as a basic template for the description of a
general approach to such problems, in particular, cases involving
zero curvature. As in the classical constant-density lattice-point
problem, there are special instances of the positive curvature case,
e.g., arithmetically defined positive definite quadratic forms, in
which the general estimate can be improved by exploiting the
underlying arithmetic character of the associated surface. In the
case in which $\partial D$ contains subsets on which the curvature
vanishes, the situation becomes vastly more intricate, although
there are general results, e.g., along the lines of the first of our
three examples in the zero curvature case. There are also arithmetic
instances of the zero curvature case having special features, as in
the second of our three examples. One can give general results for
polyhedra as well, e.g., along the lines of \cite{randol97},
\cite{skriganov98}, \cite{tarnopol2}. The methods will in general
mimic those for the constant-density case, once one is in possession
of the appropriate Fourier transform asymptotics (for the
asymptotics in the polyhedral case, cf.\ \cite{randol84b}). A kind
of meta-conclusion is that in general, the derivable asymptotics
associated with the presently considered class of problems coincide
with the corresponding results for the classical constant-density
case, by virtue of the fact that the relevant Fourier transform
asymptotics are effectively identical.

\newcommand{\noopsort}[1]{}

\bigskip

\begin{flushleft}
{\sc Ph.D Program in Mathematics\\CUNY Graduate Center\\365 Fifth
Avenue\\New York, NY 10016}
\end{flushleft}


\begin{thebibliography}{10}

\bibitem{brandolini1}
L.~Brandolini, L.~Colzani, and A.~Torlaschi.
\newblock Mean square decay of {F}ourier transforms in {E}uclidean and non
  {E}uclidean spaces.
\newblock {\em Tohuku Math. J.}, 53(3):467--478, 2001.

\bibitem{zelditch1}
Michael~R. Douglas, Bernard Shiffman, and Steve Zelditch.
\newblock Critical points and supersymmetric vacua, iii: string/m models.
\newblock {\em arXiv:math-ph/0506015}, 1, 2005.

\bibitem{hlawka1}
Edmund Hlawka.
\newblock {\"{U}}ber {I}ntegrale auf konvexen {K\"{o}}rpern {I}.
\newblock {\em Monatsh. Math.}, 54:1--36, 1950.

\bibitem{khintchine}
A.~Khintchine.
\newblock Ein {S}atz \"uber {K}ettenbr\"uche, mit arithmetischen {A}nwendungen.
\newblock {\em Math. Zeitschrift}, 18:289--306, 1923.

\bibitem{randol66a}
Burton Randol.
\newblock A lattice-point problem.
\newblock {\em Trans. Amer. Math. Soc.}, 121(1):257--268, 1966.

\bibitem{randol66b}
Burton Randol.
\newblock A lattice-point problem {II}.
\newblock {\em Trans. Amer. Math. Soc.}, 125:101--113, 1966.

\bibitem{randol69a}
Burton Randol.
\newblock On the {F}ourier transform of the indicator function of a planar set.
\newblock {\em Trans. Amer. Math. Soc.}, 139:271--276, {\noopsort{1969a}}1969.

\bibitem{randol69b}
Burton Randol.
\newblock On the asymptotic behavior of the {F}ourier transform of a convex
  set.
\newblock {\em Trans. Amer. Math. Soc.}, 139:279--285, {\noopsort{1969b}}1969.

\bibitem{randol84b}
Burton Randol.
\newblock The behavior under projection of dilating sets in a covering space.
\newblock {\em Trans. Amer. Math. Soc.}, 285:855--859, 1984.

\bibitem{randol97}
Burton Randol.
\newblock On the number of integral lattice-points in dilations of polyhedra.
\newblock {\em Int. Math. Res. Not.}, (6):259--270, 1997.

\bibitem{skriganov93}
Maxim Skriganov.
\newblock On integer points in polygons.
\newblock {\em Annales de l'{I}nstitut {F}ourier}, 43:313--323, 1993.

\bibitem{skriganov98}
Maxim Skriganov.
\newblock Ergodic theorems on \mbox{$SL(n)$}, {D}iophantine approximations and
  anomalies in the lattice-point problem.
\newblock {\em Invent. Math.}, 132(1):1--72, 1998.

\bibitem{skriganov00}
Maxim Skriganov and A.N. Starkov.
\newblock On logarithmically small errors in the lattice point problem.
\newblock {\em Ergodic Theory Dynam. Systems}, 20(5):1469--1476, 2000.

\bibitem{svensson}
Ingvar Svensson.
\newblock Estimates for the {F}ourier transform of the characteristic function
  of a convex set.
\newblock {\em Ark. Mat.}, 9:11--22, 1971.

\bibitem{varchenko1}
A.N. Varchenko.
\newblock On the number of lattice points in a domain.
\newblock {\em Uspekhi Mat. Nauk}, 37(3):177--178, 1982.

\bibitem{tarnopol2}
Marysia~Tarnopolska Weiss.
\newblock On the number of lattice-points in a compact $n$-dimensional
  polyhedron.
\newblock {\em Proc. Amer. Math. Soc.}, 74(1):124--127, 1979.

\end{thebibliography}
\end{document}